\documentclass[11pt]{amsart}

\usepackage{amsmath,amsthm, amssymb}
\usepackage[T1]{fontenc}
\usepackage{latexsym}
\usepackage{subfigure, epsfig, epsf}
\usepackage{color,fancyhdr,lastpage,graphicx}
\usepackage{xspace}
\usepackage{wasysym}
\usepackage{setspace}
\usepackage{bbm}
\usepackage{stmaryrd}

\usepackage{pstricks}
\usepackage{pst-node}
\usepackage{pst-tree}

\newcommand{\RR}{\mathbb{R}}

 \newcommand{\mcC}{\mathcal{C}}   
    
        \newcommand{\mcG}{\mathcal{G}}

  \newcommand{\be}{{\bf e}}

    \newcommand{\bfX}{{\bf X}}  \newcommand{\bfY}{\bf Y}

\newcommand{\st}{such that }

\newcommand{\ep}{\epsilon}
\renewcommand{\leq}{\leqslant}
\renewcommand{\geq}{\geqslant}

\renewcommand{\st}{such that }

\newcommand{\ssk}{\smallskip}

\newtheorem{thm}{\hspace{-0.15cm}  {\sc Theorem} }
\newtheorem{prop}[thm]{\hspace{-0.15cm} {\sc Proposition}}
\newtheorem{defn}[thm]{ \hspace{-0.3cm} {\sc Definition}}

\numberwithin{equation}{section} 

\newenvironment{Dem}{%
    \begin{list}{\hspace{0.5cm}{\sc Proof --}}{%
        \setlength{\topsep}{0pt}%
        \setlength{\leftmargin}{0pt}%
        \setlength{\rightmargin}{0pt}%
        \setlength{\listparindent}{0pt}%
        \setlength{\itemindent}{0pt}%
        \setlength{\parsep}{0pt}%
        \addtolength{\leftmargin}{20pt}%
        \addtolength{\rightmargin}{0pt}%
    } \item }{\hfill{\space $\rhd$}\end{list}\smallskip}

    \newenvironment{DemPropEstimatesGeneralRDE}{%
    \begin{list}{\hspace{0.5cm}{\sc Proof of proposition \ref{PropFundamentalEstimate} --}}{%
        \setlength{\topsep}{0pt}%
        \setlength{\leftmargin}{0pt}%
        \setlength{\rightmargin}{0pt}%
        \setlength{\listparindent}{0pt}%
        \setlength{\itemindent}{0pt}%
        \setlength{\parsep}{0pt}%
        \addtolength{\leftmargin}{20pt}%
        \addtolength{\rightmargin}{0pt}%
    } \item }{\hfill{\space $\rhd$}\end{list}\smallskip}

\title{Flows driven by Banach space-valued rough paths}
\date{\today}
\author{I. Bailleul}
\address{IRMAR, 263 Avenue du General Leclerc, 35042 RENNES, France}
\email{ismael.bailleul@univ-rennes1.fr}
\thanks{This research was partially supported by an ANR grant "Retour post-doctorant".}

\begin{document}

\maketitle

\begin{abstract}
We show in this note how the machinery of $\mcC^1$-approximate flows devised in the work {\it Flows driven by rough paths}, and applied there to reprove and extend most of the results on Banach space-valued rough differential equations driven by a finite dimensional rough path can be used to deal with rough differential equations driven by an infinite dimensional Banach space-valued weak geometric H\"older $p$-rough paths, for any $p>2$, giving back Lyons'  theory in its full force in a simple way.
\end{abstract}

\section{Introduction}
\label{SectionIntroduction}

\subsection{Rough differential equations}
\label{SubsectionIntroRDE}

The theory of rough paths invented in the mid 90' by Terry Lyons \cite{Lyons97} has had several reformulations ever since, which have both widen its scope and made it more accessible.  Of central importance in Lyons' theory is the fact that rough signals driving systems are not accurately described by the data of their increments; one needs information on what happens to the path between any two sampling times to understand the output of the system as a function of the driving signal. The nature of this additional information makes it necessary to introduce a rich stratified algebraic structure. Rough paths are paths with values in that structure, the increments of the initial rough signal being recorded in the first level of the associated rough path. A Riemann-type integral, with any rough path as integrator, can be defined in that setting, which enables to formulate differential equations driven by rough paths as fixed point problems to integral equations in the enriched algebraic structure. Gubinelli extracted in \cite{Gubinelli} the core of Lyons' machinery in the notion of controlled path. He put forward the essential fact that the class of integrands for which Riemann-type integrals can be defined in the context of rough paths is made up of time-dependent maps whose time increments are time-dependent linear transforms of the first level of the driving rough path, up to some smoother term. Gubinelli's framework was used in several works (e.g. \cite{DeyaGubinelliTindel, GubinelliKdV, FlandoliGubinelliPriola, GubinelliTindel}) to deal with a variety of subtle and hard problems in stochastic partial differential equations, and was one of the seeds of Hairer's impressive works \cite{HairerKPZ, HairerRegularityStructures}, where he lays the foundations of a new theoretical framework for the study of stochastic partial differential equations and solves some longstanding fundamental problems.

\medskip

In parallel to this line of development, Friz and Victoir developped in \cite{FVEuler} and \cite{FVBook} a purely dynamical approach to rough differential equations, based on Davie's insighful work \cite{Davie}. Let $\textrm{F}=(V_1,\dots,V_\ell)$ be a finite collection of smooth enough vector fields on $\RR^d$ and $\bfX$ be a rough path over $\RR^\ell$. In Friz-Victoir's approach, a solution to the rough differential equation
$$
dx_u = \textrm{F}(x_u)d{\bf X}_u
$$
on $\RR^d$, is the limit in uniform topology of the solutions to ordinary differential equations of the form
\begin{equation}
\label{EqFVODE}
\dot y_u = \sum_{i=1}^\ell V_i(y_u)\dot h^i_u,
\end{equation}
for some smooth control $h$, when the rough path canonically associated to $h$ converges in a rough path sense to $\bfX$. The technical core of their approach takes the form of some estimates on some high order Euler expansion to the solution to equation \eqref{EqFVODE} involving only some $p$-variation norm of $h$. However, both \cite{Davie} and \cite{FVEuler} obtain dimension-dependent estimates which prevent a direct extension of their results to investigate rough differential equations driven by infinite dimensional rough paths, as can be done in Lyons or Gubinelli's framework.

\medskip

A different approach to rough differential equations was proposed recently in the work "\textit{Flows driven by rough paths}", \cite{RPDrivenFlows}. Given some finite collection $\textrm{F}=(V_1,\dots,V_\ell)$ of smooth enough vector fields on a Banach space $U$, a flow $\big(\varphi_{ts}\big)_{0\leq s\leq t\leq T}$ of maps from $U$ to itself is said to be a solution of the equation
\begin{equation}
\label{EqFlowRDE}
d\varphi_s = \textrm{F}{\bfX}(ds)
\end{equation}
if its increments can accurately be described by some awaited Euler-type expansion, as in Davie or Friz-Victoir's approach. The introduction of the notion of $\mcC^1$-approximate flow, and the proof that a $\mcC^1$-approximate flow determines a unique flow close enough to it, made it possible to recover and extend most of the basic results on Banach space valued rough differential equations driven by a finite dimensional rough signal, under optimal regularity assumtions on the vector fields $V_i$. Existence, well-posedness, Taylor expansion and non-explosion under linear growth conditions of the vector fields were proved in \cite{RPDrivenFlows}.

\ssk

We show in the present note that the machinery of $\mcC^1$-approximate flows applies equally well to prove well-posedness results for Banach space valued rough differential equations driven by infinite dimension rough signals, extending to that setting Davie-Friz-Victoir's Euler estimates as a consequence. This infinite dimensional version of Davie's generalized estimate was proved recently in \cite{BoutaibGyurkoLyonsYang} by adapting the method of proof of Davie \cite{Davie} and Friz-Victoir \cite{FVEuler}, and using Lyons' universal limit theorem. Our method yields both Lyons' theorem and Davie's estimate at a time.

\medskip

We recall in section \ref{SubsectionApproximateFlows} how to construct flows from $\mcC^1$-approximate flows \cite{RPDrivenFlows}, and construct in section \ref{SectionInfiniteDimRPDrivenFlows} a $\mcC^1$-approximate flow with the awaited Euler expansion. Well-posedness of the rough differential equation \eqref{EqFlowRDE} on flows and Davie's estimate are direct consequences of the result on $\mcC^1$-approximate flows recalled in section \ref{SubsectionApproximateFlows}.

\medskip

The construction of a solution flow to equation \eqref{EqFlowRDE} done in section \ref{SectionInfiniteDimRPDrivenFlows} requires three ingredients, introduced in subsections \ref{SubsectionApproximateFlows}, \ref{SubsectionInfiniteDimRP} and \ref{SubsectionDifferentialOperators}. Throughout that work, we denote $U$ and $V$ two Banach spaces. Given a non-negative real number $\gamma$, we denote by $\mcC^\gamma(U,U)$ the set of $\gamma$-Lipschitz maps from $U$ to itself, in the sense of Stein, equipped with its natural norm. (See \cite{Lyons97} or \cite{LyonsStFlour}.) The notation $\mcC^\gamma(U)$ stands for the setof $\gamma$-Lipschitz real valued functions on $U$. We use the symbol $c$ to denote a constant whose value is unimportant and may change from place to place.

\medskip

\subsection{Flows and approximate flows}
\label{SubsectionApproximateFlows}

Recall that a \textit{flow} on $U$ is a family $\big(\varphi_{ts}\big)_{0\leq s\leq t\leq T}$ of maps from $U$ to itself such that $\varphi_{ss}=\textrm{Id}$, for all $0\leq s\leq T$, and $\varphi_{tu}\circ\varphi_{us} = \varphi_{ts}$, for all $0\leq s\leq u\leq t\leq T$. The notion of $\mcC^1$-approximate flow introduced in \cite{RPDrivenFlows} provides a convenient tool for constructing flows from families of maps which are almost flows. It is an elaboration upon Lyons' almost-multiplicative functionals, as understood by Gubinelli \cite{Gubinelli} or Feyel-de la Pradelle \cite{FdlP1,FdlP2}, with their sewing lemmas. We state it here under a simple form which will be enough to illustrate our purpose in section \ref{SectionInfiniteDimRPDrivenFlows}.

\begin{defn}
\label{DefnC1ApproximateFlow}
A \emph{{\bf $\mcC^1$-approximate flow}} is a family $\big(\mu_{ts}\big)_{0\leq s\leq t\leq T}$ of $\mcC^2$ maps from $U$ to itself, depending continuously on $(s,t)$ in the topology of uniform convergence, such that 
   \begin{equation}
   \label{EqRegularityBounds}
   \big\|\mu_{ts} - \textrm{\emph{Id}}\big\|_{\mcC^2} = o_{t-s}(1),
   \end{equation} 
and there exists some positive constants $c_1$ and $a>1$, such that the inequality
\begin{equation}
\label{EqMuMu}
\big\|\mu_{tu}\circ\mu_{us}-\mu_{ts}\big\|_{\mcC^1} \leq c_1\,|t-s|^a
\end{equation}
holds for any $0\leq s\leq u\leq t\leq T$.
\end{defn}

Given a partition $\pi_{ts}=\{s=s_0<s_1<\cdots<s_{n-1}<s_n=t\}$ of $(s,t)\subset [0,T]$, set 
$$
\mu_{\pi_{ts}} = \mu_{s_ns_{n-1}}\circ\cdots\circ\mu_{s_1s_0}.
$$
The following statement, proved in \cite{RPDrivenFlows}, provides a flexible tool for constructing flows from $\mcC^1$-approximate flows, while generalizing the above mentionned sewing lemmas to the non-commutative setting of $\mcC^2$ maps on $U$.
 
\begin{thm}[Constructing flows on a Banach space]
\label{ThmConstructingFlows}
A $\mcC^1$-approximate flow defines a unique flow $(\varphi_{ts})_{0\leq s\leq t\leq T}$ on $U$ such that the inequality
$$
\big\|\varphi_{ts}-\mu_{ts}\big\|_\infty \leq c|t-s|^a
$$
holds for some positive constant $c$, for all $0\leq s\leq t\leq T$ sufficiently close, say $|t-s|\leq \delta$; this flow satisfies the inequality 
\begin{equation}
\label{EqApproxVarphiMu}
\big\|\varphi_{ts}-\mu_{\pi_{ts}}\big\|_\infty\leq c_1^2 \,T\,\big|\pi_{ts}\big|^{a-1}
\end{equation} 
for any partition $\pi_{ts}$ of any interval $(s,t)\subset [0,T]$, with mesh $\big|\pi_{ts}\big|\leq \delta$. 
\end{thm}

A slightly more elaborated form of this theorem was used in \cite{RPDrivenFlows} to recover and extend most of the basic results on Banach-space valued rough differential equations driven by a finite dimensional rough signal, under optimal regularity condition on the driving vector fields. Existence, well-posedness, high order Euler expansion and non-explosion under linear growth conditions on the driving vector fields were proved. It was also used in \cite{PathDpdtRDEs} to extend these results to the setting of rough differential equations driven by path-dependent vector fields.

\subsection{Infinite dimensional rough paths}
\label{SubsectionInfiniteDimRP}

We recall briefly in this section the setting of infinite dimensional rough paths, refering to \cite{Lyons97, LyonsQian} for more material on this subject. 

\ssk

For $n\geq 1$, let $T^n_a(V)$ stand for the truncated algebraic tensor algebra, isomorphic to $\bigoplus_{i=0}^n V^{\otimes_a i}$, where $V^{\otimes_a i}$ stands for the algebraic tensor product of $V$ with itself $i$ times, and $\RR^{\otimes 0}$ is identified to $\RR$. Multiplication of two elements $\be$ and $\be'$ of $T^n(V)$ is simply denoted by $\be\be'$. There exists norms on $T^n_a(V)$ which satisfy the inequality
$$
\|\be\be'\|\leq \|\be\|\|\be'\|, 
$$
for any pair $(\be,\be')$ of elements of $T^n(V)$, with $\|\be\|=|\be|$, if $\be\in V\subset T^n(V)$. Select any such norm and denote by $T^n(V)$ the completion of $T^n_a(V)$ with respect to that norm. We have $T^n(V) = \bigoplus_{i=0}^n V^{\otimes i}$, where $V^{\otimes i}$ stands for the completion of $V^{\otimes_a i}$ with respect to the restriction to $V^{\otimes_a i}$ of the norm $\|\cdot\|$, and $V^{\otimes 1} = V$, since $V$ is complete.

Given $0\leq k\leq n$, denote by $\pi_k$ the projection map from $T^n(V)$ to $V^{\otimes k}$, and set $T^n_1(V) =\{{\be}\in T^n(V)\,;\,\pi_0({\be}) = 1\}$, with a similar definition of $T^n_0(V)$. The exponential map is a polynomial diffeomorphism from $T^n_0(V)$ to $T^n_1(V)$, with inverse map the logarithm. Let $\frak{g}^n(V)$ stand for the Lie algebra in $T^n(V)$ generated by $V$, and write $\mcG^n(V)$ for $\exp\big(\frak{g}^n(V)\big)$. 

\ssk

In this setting, and for $2\leq p$, a {\bf $V$-valued weak geometric H\"older $p$-rough path} is a $\mcG^{[p]}(V)$-valued path $\bfX$, parametrized by some time interval $[0,T]$, \st 
\begin{equation}
\label{ConditionsHolder}
\underset{0\leq s<t\leq T}{\sup}\,\frac{\big\|\pi_i {\bf X}_{ts}\big\|}{|t-s|^{\frac{i}{p}}}<\infty,
 \end{equation}
for all $i=1\dots [p]$, where we set ${\bf X}_{ts} = {\bf X}_s^{-1}{\bf X}_t$. We define the norm of $\bfX$ to be 
\begin{equation}
\llparenthesis{\bfX}\rrparenthesis = \underset{i=1\dots [p]}{\max}\;\underset{0\leq s<t\leq T}{\sup}\,\frac{\big\|\pi_i {\bf X}_{ts}\big\|}{|t-s|^{\frac{i}{p}}},
\end{equation}
and a distance $d(\bf X,\bf Y) = \llparenthesis\bfX-\bfY\rrparenthesis$, on the set of $V$-valued weak geometric H\"older $p$-rough path. This definition depends on the choice of tensor norm made above.

\subsection{Differential operators}
\label{SubsectionDifferentialOperators}

Let $\gamma>0$ be given, and F be a continuous linear map from $V$ to $\mcC^\gamma(U,U)$. For any $v\in V$, we identify the $\mcC^\gamma$ vector field $\textrm{F}(v)$ on $U$ with the first order differential operator
$$
g\in \mcC^1(U)\mapsto (D_{\cdot} g)\big(\textrm{F}(v)(\cdot)\big)\in\mcC^0(U);
$$
in those terms, we recover the vector field as $\textrm{F}(v)\textrm{Id}$. The map F is extended to the algebraic setting of $T^{[\gamma]}(V)$ by setting $\textrm{F}^\otimes(1) = \textrm{Id} : \mcC^0(U)\mapsto \mcC^0(U)$, and defining $\textrm{F}^\otimes (v_1\otimes\cdots\otimes v_k)$, for all $1\leq k\leq [\gamma]$ and $v_1\otimes\cdots\otimes v_k\in V^{\otimes k}$, as the $k^{\textrm{th}}$-order differential operator from $\mcC^k(U)$ to $\mcC^0(U)$, defined by the formula
$$
\textrm{F}^\otimes (v_1\otimes\cdots\otimes v_k) = \textrm{F}\big(v_1\big)\cdots \textrm{F}\big(v_k\big),
$$
and by requiring linearity. The choice of a tensor norm on $T^{[\gamma]}(V)$ ensures that this definition extends continuously to the restriction of the completed space $V^{\otimes k}$ as a continuous linear map from $V^{\otimes k}$ to the set of $k^\textrm{th}$-order differential operators from $\mcC^k(U)$ to $\mcC^0(U)$, endowed with the operator norm. So, for $\be\in\frak{g}^n(V)$ the differential operator $\textrm{F}(\be)$ is actually a first order differential operator, and we have
\begin{equation}
\label{EqMorphismF}
\textrm{F}^\otimes(\be)\,\textrm{F}^\otimes(\be') = \textrm{F}^\otimes(\be\be'),
\end{equation}
for all $\be,\be'\in T^n(V)$.

\section{Flows driven by Banach space-valued rough paths}
\label{SectionInfiniteDimRPDrivenFlows}

We show in this section how to associate a $\mcC^1$-approximate flow to some $\mcG^{[p]}(V)$-valued weak geometric H\"older $p$-rough path and some $V$-dependent vector field F as in section \ref{SubsectionDifferentialOperators}, in such a way that it has the awaited Euler expansion. A solution flow to the equation 
\begin{equation}
\label{EqRDE}
d\varphi_s = \textrm{F}{\bfX}(ds)
\end{equation}
will be defined as a flow uniformly close to the $\mcC^1$-approximate flow. Well-posedness of the rough differential equation \eqref{EqRDE} on flows will follow as a direct consequence of theorem \ref{ThmConstructingFlows}. 

\bigskip

Let $2\leq p$ be given, together with a $\mcG^{[p]}(V)$-valued weak-geometric H\"older $p$-rough path $\bfX$, defined on some time interval $[0,T]$, and some continuous linear map F from $V$ to $\mcC^{[p]+1}(U,U)$. For any $0\leq s\leq t\leq T$, denote by ${\bf \Lambda}_{ts}$ the logarithm of ${\bfX}_{ts}$,  and let $\mu_{ts}$ stand for the well-defined time $1$ map associated with the ordinary differential equation
\begin{equation}
\label{EqRPODE}
\dot y_u = \textrm{F}^\otimes\big({\bf\Lambda}_{ts}\big)(y_u),\quad 0\leq u\leq 1.
\end{equation}
This equation is indeed an ordinary differential equation as ${\bf\Lambda}_{ts}$ is an element of $\frak{g}^{[p]}(V)$. It is a consequence of classical results from ordinary differential equations, and the definition of the norm on the space of weak-geometric H\"older $p$-rough paths, that the solution map $(r,x)\mapsto y_r$, with $y_0=x$, depends continuously on $\big((s,t),\bfX\big)$, and satisfies the following basic estimate. We have

\begin{equation}
\label{EqRoughODEsEstimate}
\big\|y_r-\textrm{Id}\big\|_{\mcC^1} \leq c\left(1+\llparenthesis{\bfX}\rrparenthesis^{[p]}\right)|t-s|^\frac{1}{p}, \quad 0\leq r\leq 1.
\end{equation}
The following elementary proposition refines part of the above estimate, and is our basic step for studying flows driven by Banach space-valued weak geometric H\"older $p$-rough paths. 

\begin{prop}
\label{PropFundamentalEstimate}
There exists a positive constant $c$, depending only on the data of the problem, \st the inequality
\begin{equation}
\label{EqFundamentalEstimate}
 \Big\|f\circ\mu_{ts} - \textrm{\emph{F}}^\otimes\big({\bfX}_{ts}\big)f\Big\|_\infty \leq c\Big(1+\llparenthesis {\bfX} \rrparenthesis^{[p]}\Big)\,\|f\|_{[p]+1}\,|t-s|^{\frac{[p]+1}{p}}
\end{equation}
holds for any $f\in\mcC^{[p]+1}(U)$. 
\end{prop}

This proposition was first proved under this form in \cite{RPDrivenFlows}, following the same pattern of proof as below, for rough differential equations driven by a finite dimensional rough path. The above infinite dimensional version of it, with $f=\textrm{Id}$, was proved independently in \cite{BoutaibGyurkoLyonsYang} using the same method.

\ssk

The proof of this proposition and the following one are based on the elementary identity \eqref{EqExactFormulaMuTs} below, obtained by applying repeatedly the identity
$$
f\big(y_r\big) = f(x) + \int_0^r\Big(\textrm{F}^\otimes\big({\bf\Lambda}_{ts}\big)f\Big)(y_u)\,du,\quad 0\leq r\leq 1,
$$
together with the morphism property \eqref{EqMorphismF}, and by separating the terms according to their size in $|t-s|$. Set $\Delta_n := \big\{(s_1,\dots,s_n)\in [0,T]^n\,;\,s_1\leq \cdots\leq s_n\big\}$, for $2\leq n\leq [p]$. The summation below is above indices $k_i$ with range in $\big\llbracket 1,[p]\big\rrbracket$; we also write $ds$ for $ds_n\dots ds_1$. For a $\gamma$-Lipschitz function $f$,  we have
\begin{equation}
\label{EqExactFormulaMuTs}
\begin{split}
&f\big(\mu_{ts}(x)\big) = f(x) + \sum_{\ell=1}^n \frac{1}{\ell !} \sum_{k_1+\cdots+k_\ell \leq [p]} \Big(\textrm{F}^\otimes\big(\pi_{k_\ell}{\bf\Lambda}_{ts}\cdots \pi_{k_1}{\bf\Lambda}_{ts}\big)f\Big)(x) \\
&+\sum_{k_1+\cdots+k_n \leq [p]} \int_{\Delta_n} \Big\{\Big(\textrm{F}^\otimes\big(\pi_{k_n}{\bf\Lambda}_{ts}\cdots\pi_{k_1}{\bf\Lambda}_{ts}\big)f\Big)\big(y_{s_n}\big) -  \Big(\textrm{F}^\otimes \big(\pi_{k_n}{\bf\Lambda}_{ts}\cdots\pi_{k_1}{\bf\Lambda}_{ts}\big)f\Big)(x)\Big\}\,ds \\
&+ \sum_{\ell=1}^n \frac{1}{\ell !} \sum_{k_1+\cdots+k_\ell \geq [p]+1} \Big(\textrm{F}^\otimes\big(\pi_{k_\ell}{\bf\Lambda}_{ts}\cdots \pi_{k_1}{\bf\Lambda}_{ts}\big)f\Big)(x) \\
&+\sum_{k_1+\cdots+k_n \geq [p]+1} \int_{\Delta_n} \Big\{\Big(\textrm{F}^\otimes\big(\pi_{k_n}{\bf\Lambda}_{ts}\cdots \pi_{k_1}{\bf\Lambda}_{ts}\big)f\Big)\big(y_{s_n}\big) -  \Big(\textrm{F}^\otimes\big(\pi_{k_n}{\bf\Lambda}_{ts}\cdots \pi_{k_1}{\bf\Lambda}_{ts}\big)f\Big)(x)\Big\}\,ds.
\end{split}
\end{equation}
\noindent We denote by $\ep^{f\,;\,n}_{ts}(x)$ the sum of the last two lines, made up of terms of size at least $|t-s|^{[p]+1}$. In the case where $n=[p]$, the terms in the second line involve only indices $k_j$ with $k_j=1$, so the elementary estimate \eqref{EqRoughODEsEstimate} can be used to control the increment in the integral, showing that this second line is of order $|t-s|^\frac{[p]+1}{p}$, in infinite norm, as the maps $\textrm{F}^\otimes\big(\pi_{k_n}{\bf\Lambda}_{ts}\cdots\pi_{k_1}{\bf\Lambda}_{ts}\big)f$ are $\mcC^1_b$; we include it in the remainder $\ep^{f\,;\,[p]}_{ts}(x)$.

\medskip

\begin{DemPropEstimatesGeneralRDE}
Applying the above formula for $n=[p]$, together with the fact that $\exp\big({\bf \Lambda}_{ts}\big) = {\bfX}_{ts}$, we get the identity
\begin{equation*}
f\big(\mu_{ts}(x)\big) = \Big(\textrm{F}^\otimes\big({\bfX}_{ts}\big)f\Big)(x) + \ep^{f\,;\,[p]}_{ts}(x).
\end{equation*}
It is clear on the formula for $\ep^{f\,;\,[p]}_{ts}(x)$ that its absolute value is bounded above by a constant multiple of $\Big(1+\llparenthesis {\bfX} \rrparenthesis^{[p]}\Big)|t-s|^\frac{[p]+1}{p}$, for a constant depending only on the data of the problem and $f$ as in \eqref{EqFundamentalEstimate}. 
\end{DemPropEstimatesGeneralRDE}

\medskip

A further look at formula \eqref{EqExactFormulaMuTs} makes it clear that if $2\leq n\leq [p]$ and $f$ is $\mcC^{n+1}_b$, the estimate 
\begin{equation}
\label{EqEstimateC1SizeEpsilon}
\Big\|\ep^{f\,;\,n}_{ts}\Big\|_{\mcC^1} \leq c\Big(1+\llparenthesis {\bfX} \rrparenthesis^{[p]}\Big)\,\|f\|_{\mcC^{n+1}}|t-s|^\frac{[p]+1}{p},
\end{equation}
holds as a consequence of formula \eqref{EqRoughODEsEstimate}, for a constant $c$ depending only on the data of the problem. 

\begin{prop}
\label{PropSummaryPropertiesGeneralCase}
The family of maps $\big(\mu_{ts}\big)_{0\leq s\leq t\leq T}$ is a $\mcC^1$-approximate flow. 
\end{prop}

\medskip

\begin{Dem}
As F is $\mcC^{[p]+1}$ as a function on $U$, with $[p]+1\geq 3$, the inequality $\big\|\mu_{ts}-\textrm{Id}\big\|_{\mcC^2} = o_{t-s}(1)$, is given by classical results on ordinary differential equations; we turn to proving the $\mcC^1$-approximate flow property \eqref{EqMuMu}. of $\big(\mu_{ts}\big)_{0\leq s\leq t\leq T}$. We first use for that purpose formula \eqref{EqExactFormulaMuTs} to write
\begin{equation}
\begin{split}
\mu_{tu}\big(\mu_{us}(x)\big) &= \Big(\textrm{F}^\otimes\big({\bfX}_{tu}\big)\textrm{Id}\Big)\big(\mu_{us}(x)\big) +\ep^{\textrm{Id}\,;\,[p]}_{tu}\big(\mu_{us}(x)\big) \\
&= \mu_{us}(x) + \sum_{m=1}^{[p]} \Big(\textrm{F}^\otimes\big(\pi_m{\bfX}_{tu}\big)\textrm{Id}\Big)\big(\mu_{us}(x)\big) +\ep^{\textrm{Id}\,;\,[p]}_{tu}\big(\mu_{us}(x)\big)
\end{split}
\end{equation}
The remainder $\ep^{\textrm{Id}\,;\,[p]}_{tu}\big(\mu_{us}(x)\big)$ has a $\mcC^1$-norm bounded above by $c\Big(1+\llparenthesis {\bfX} \rrparenthesis^{[p]}\Big)^2|t-u|^\frac{[p]+1}{p}$, due \eqref{EqEstimateC1SizeEpsilon} and inequality \eqref{EqRoughODEsEstimate}.

\ssk

To deal with the term $\Big(\textrm{F}^\otimes\big(\pi_m{\bfX}_{ts}\big)\textrm{Id}\Big)\big(\mu_{us}(x)\big)$, we use formula \eqref{EqExactFormulaMuTs} with $n=[p]-m$. Writing $ds$ for $ds_{[p]-m}\dots ds_1$, it folllows that
\begin{equation}
\label{EqControlFOtimes}
\begin{split}
&\Big(\textrm{F}^\otimes\big(\pi_m{\bfX}_{tu}\big)\textrm{Id}\Big)\big(\mu_{us}(x)\big) = \Big(\textrm{F}^\otimes\big(\pi_m{\bfX}_{ts}\big)\textrm{Id}\Big)(x) + \sum_{\ell=1}^{[p]-m} \frac{1}{\ell !}\sum_{k_1+\cdots+k_\ell\leq [p]} \Big(\textrm{F}^\otimes\big(\pi_{k_\ell}{\bf \Lambda}_{ts}\cdots\pi_{k_1}{\bf \Lambda}_{ts}\,\pi_m{\bfX}_{tu}\big) \textrm{Id}\Big) (x)  \\ 
&+ \ep^{\bullet\,;\,{p}-m}_{us}(x) \\
&+ \sum \int_{\Delta_{[p]-m}} \Big\{\Big(\textrm{F}^\otimes\big(\pi_{k_{[p]-m}}{\bf \Lambda}_{us}\cdots\pi_{k_1}{\bf \Lambda}_{us}\,\pi_m{\bfX}_{tu}\big)\textrm{Id}\Big)\big(y_{s_{[p]-m}}\big) - \Big(\textrm{F}^\otimes\big(\pi_{k_\ell}{\bf \Lambda}_{us}\cdots\pi_{k_1}{\bf \Lambda}_{us}\,\pi_m{\bfX}_{tu}\big)\textrm{Id}\Big)(x)\Big\}\,ds,
\end{split}
\end{equation}
where the last sum is over the set $\big\{k_1+\dots+k_{[p]-m}\leq [p]\big\}$ of indices. The notation $\bullet$ in the above identity stands for the $\mcC^{[p]+2-m}_b$ function $\textrm{F}^\otimes\big(\pi_m{\bfX}_{tu}\big)\textrm{Id}$ on $U$; it has $\mcC^1$-norm controlled by \eqref{EqEstimateC1SizeEpsilon}. It is then straighforward to use \eqref{EqRoughODEsEstimate} to bound above the $\mcC^1$-norm of the third line in equation \eqref{EqControlFOtimes} by a constant multiple of $\Big(1+\llparenthesis {\bfX} \rrparenthesis^{[p]}\Big)\,|t-s|^\frac{[p]+1}{p}$. Writing 
$$
\mu_{us}(x) = \Big(\textrm{F}^\otimes\big({\bfX}_{us}\big)\textrm{Id}\Big)(x) + \ep^{\textrm{Id}\,;\,[p]}_{us}(x),
$$
and using the identities $\exp\big({\bf \Lambda}_{us}\big)=\bfX_{us}$ and ${\bfX}_{ts} = {\bfX}_{us}{\bfX}_{tu}$, we see that
$$
\mu_{tu}\big(\mu_{us}(x)\big) = \mu_{ts}(x) + \ep_{ts}(x),
$$
with a remainder
\begin{equation}
\label{EqEstimateRemainder} 
\|\ep_{ts}\|_{\mcC^1} \leq c\Big(1+\llparenthesis {\bfX} \rrparenthesis^{[p]}\Big)|t-s|^\frac{[p]+1}{p}.
\end{equation}
\end{Dem}

\medskip

In view of proposition \ref{PropFundamentalEstimate}, the following definition is to be thought of as an analogue of Davie's definition \cite{Davie} of a solution to a classical rough differential equation, in terms of Euler expansion.

\begin{defn}
\label{DefnGeneralRDESolution}
A \textbf{\emph{flow}} $(\varphi_{ts}\,;\,0\leq s\leq t\leq T)$ is said to \textbf{\emph{solve the rough differential equation}}
\begin{equation}
\label{RDEGeneral}
d\varphi = \textrm{\emph{F}}\,\bfX(dt)
\end{equation}
if there exists a constant $a>1$ independent of $\bfX$ and two possibly $\bfX$-dependent positive constants $\delta$ and $c$ such that
\begin{equation}
\label{DefnSolRDEGeneral}
\|\varphi_{ts}-\mu_{ts}\|_\infty \leq c\,|t-s|^a
\end{equation}
holds for all $0\leq s\leq t\leq T$ with $t-s\leq\delta$.
\end{defn}

It is clear on this definition that if $x\in U$ and $\big(\varphi_{ts}\big)_{0\leq s\leq t\leq T}$ is a solution flow to equation \eqref{RDEGeneral}, then the trajectory $\big(\varphi_{t0}(x)\big)_{0\leq t\leq T}$ is the first level of a solution to the classical rough differential equation $dx_t = \textrm{F}d{\bfX}_t$, in Lyons' sense. The following well-posedness result follows directly from theorem \ref{ThmConstructingFlows} and proposition \ref{PropSummaryPropertiesGeneralCase}.

\begin{thm}
\label{ThmMainResultGeneral}
Suppose \emph{F} is a continuous linear map from $V$ to $\mcC^{[p]+1}(U,U)$; extend it as in section \ref{SubsectionDifferentialOperators} into a differential operator. Then the rough differential equation $d\varphi = \textrm{\emph{F}}\,\bfX(dt)$ has a unique solution flow; it depends continuously on $\big((s,t),\bfX\big)$. 
\end{thm}

\medskip

\begin{Dem}
\noindent With the notations of definition \ref{DefnC1ApproximateFlow}, identity \eqref{EqEstimateRemainder} means that equation \eqref{EqMuMu} holds with $c_1=c\Big(1+\llparenthesis {\bfX} \rrparenthesis^{[p]}\Big)$ and $a=\frac{[p]+1}{p}$. So theorem \ref{ThmConstructingFlows} ensures the existence of a unique flow $\big(\varphi_{ts}\big)_{0\leq s\leq t\leq T}$ close enough to $\big(\mu_{ts}\big)_{0\leq s\leq t\leq T}$; it further satisfies the inequality
\begin{equation}
\label{EqApproxPhiMu}
\Big\|\varphi_{ts}-\mu_{\pi_{ts}}\Big\|_\infty \leq c\Big(1+\llparenthesis {\bfX} \rrparenthesis^{[p]}\Big)^2T\,\big|\pi_{ts}\big|^\frac{1}{p},
\end{equation}
for any partition $\pi_{ts}$ of $(s,t)\subset [0,T]$ of mesh $\big|\pi_{ts}\big|\leq \delta$, as a consequence of inequality \eqref{EqApproxVarphiMu}. As these bounds are uniform in $(s,t)$, and for $\bfX$ in a bounded set of the space of H\"older $p$-rough paths, and each $\mu_{\pi_{ts}}$ is a continuous function of $\big((s,t),\bfX\big)$, as a composition of continuous functions, the flow $\varphi$ depends continuously on $\big((s,t),\bfX\big)$.
\end{Dem}

\ssk

Davie-Friz-Victoir's estimate follows from proposition \ref{PropFundamentalEstimate} and theorem \ref{ThmMainResultGeneral} under the form 
$$
\big\|\varphi_{ts}-\textrm{F}^\otimes\big({\bfX}_{ts}\big)\textrm{Id}\big\|_\infty \leq c\Big(1+\llparenthesis {\bfX} \rrparenthesis^{[p]}\Big)\,|t-s|^\frac{\gamma}{p}.
$$
By using the refined definition of a $\mcC^1$-approximate flow given in \cite{RPDrivenFlows} and the analogue of theorem \ref{ThmConstructingFlows} which holds for it, the above method can be used to prove theorem \ref{ThmMainResultGeneral} under the essentially optimal condition that F takes values in $\mcC^\gamma(U,U)$, for any choice of $2<p<\gamma\leq [p]+1$. The other results proved in \cite{RPDrivenFlows}: high order Euler expansion, non-explosion under linear growth conditions on F, a Peano theorem on existence of a solution, can be proved in the present setting as well, by a straightforward adaptation of the notations of \cite{RPDrivenFlows}.

\vfill

\end{document}